\documentclass[12pt]{article}
 \usepackage{amsfonts,amssymb,amsmath,mathrsfs}
 \usepackage{color}
 \usepackage[svgnames]{xcolor}

 \newtheorem{theorem}{Theorem}[section]
 \newtheorem{lemma}[theorem]{Lemma}
 \newtheorem{corol}[theorem]{Corollary}
 \newtheorem{prop}[theorem]{Proposition}
 \newtheorem{remark}[theorem]{Remark}
 \newtheorem{example}[theorem]{Example}
 \newtheorem{condition}[theorem]{Condition}

 \def\btheorem{\begin{theorem}\sl{}
 \def\etheorem{\end{theorem}}}
 \def\blemma{\begin{lemma}\sl{}
 \def\elemma{\end{lemma}}}

 \def\bremark{\begin{remark}\sl{}
 \def\eremark{\end{remark}}}

 \def\benumerate{\begin{enumerate}}\def\eenumerate{\end{enumerate}}
 \def\bitemize{\begin{itemize}}\def\eitemize{\end{itemize}}
 \def\itm{\item}

 \def\beqlb{\begin{eqnarray}}
 \def\eeqlb{\end{eqnarray}}
 \def\beqnn{\begin{eqnarray*}}
 \def\eeqnn{\end{eqnarray*}}

 \def\pf{\noindent{\it Proof.~~}}
 \def\qed{\hfill$\Box$\medskip}

 \def\<{\langle}\def\>{\rangle}

 \def\mcr{\mathscr}\def\mbb{\mathbb}
 \def\mbf{\mathbf}\def\mrm{\mathrm}

 \def\ar{\!\!&}

 \def\d{\mrm{d}}\def\e{\mrm{e}}

 \def\lline{--------------------------------------}

\begin{document}
\bigskip\bigskip

\centerline{\LARGE\bf Scaling Limits of Controlled Branching Processes} 

\bigskip\bigskip

\centerline{ Jiawei Liu}

\bigskip

\centerline{School of Mathematical Sciences, Beijing Normal University,}

\centerline{Beijing 100875, People's Republic of China}

\centerline{E-mails:  {\tt
jwliu@mail.bnu.edu.cn}}

\bigskip

\centerline{\lline\lline\lline}

\medskip

{\narrower

\noindent\textit{Abstract.} In this paper, a special sequence of controlled
branching processes is considered. We provide a simple set of sufficient
conditions for the weak convergence of such processes to a weak solution to
a kind of continuous branching processes with dependent immigration.

\smallskip


\smallskip


\par}

\bigskip

\centerline{\lline\lline\lline}

\bigskip

\section{Introduction}

 Suppose that there is a family of random variables
 $\{\xi_{n,i}: n,i=1,2,\cdots\}$ with values
 in $\mbb{N}:=\{ 0,1,2,\cdots\}$, which are mutually independent. Given an $\mbb{N}$-valued
 random variable $Z(0)$ independent of $\{\xi_{n,i}\}$, a
 \textit{Galton-Watson process} (GW process)
 $\{Z(n):n\in \mbb{N}\}$ can be inductively defined by
 \beqlb\label{1.1}
Z(n+1)=\sum^{Z(n)}_{i=1} \xi_{n,i},\ \quad n= 0, 1, 2,\cdots.
 \eeqlb
 Here, we understand $\sum_{i=1}^0=0$. In the classical GW process, $\xi_{n,i}$ is viewed as the offspring produced by the $i$-th individual in $n$-th generation. It was proved that the scaling limits of GW processes can be a \textit{continuous-state branching process }(CB process);
 see, e.g., \cite{Ka71,Li06}. For more theories on scaling limits of generalized GW processes, one can refer to \cite{BaC19,Fa22,Pa16}.
While each individual is influenced by environments, the mechanism of reproduction may vary in different generations. The mechanism here is understood as a competition or a interaction, and it's initiated from  \cite{Se74}, where
 Sevastyanov and Zubkov generalized the model of GW processes by considering a constant control on the growth of population size at each generation. Later, Yanev \cite{Ya76} consider the conditions that the controls are random and i.i.d. He introduced a model of \textit{controlled branching process with random control function} (CBP), which can be formulated as follows. Let $\{ \phi^{(n)}(i) : i= 0, 1, 2,\cdots\}, n= 0, 1, 2,\cdots$ be a mutually independent random function having the same distribution for each $n$.  A CBP $\{Z(n):n\in\mbb{N}\}$ was constructed inductively as follows:
 \beqlb\label{1.2}
 Z(n+1)=\sum^{\phi^{(n)}(Z(n))}_{i=1} \xi_{n,i},\ \quad n= 0, 1, 2,\cdots.
 \eeqlb
 From equality above, $\phi^{(n)}(i)$ is viewed as a random control function. The probabilistic theory on this model was developed by Gonz\'alez et al. \cite{Go02,Go03}
 and so on.

 Considering a sequence of such processes $\{Z_k(n);~n\geq 0\}_{k\geq 1}$ with offspring and random control functions $\big(\xi_{n,i}^{(k)},\phi_k^{(n)}(i)\big)$, we concentrate on how the rescaled processes $\{Y_k(t):=Z_k(\lfloor \gamma_k t\rfloor)/k;~t\geq 0\}_{k\geq 1}$ converges on the Skorohod space as $k\rightarrow\infty$, where $\gamma_k$ is a sequence of positive increasing constants tending to $\infty$. Such a question was partially answered by Gonz\'alez and del Puerto
 \cite{Go12}. They proved the weak convergence of a sequence of CBPs to a diffusion
 process under some restrictions on the means and variances of $\xi_{n,i}^{(k)}$ and  $\phi_k^{(n)}(i)$. In fact, the diffusion that they obtained is a Feller branching diffusion with immigration. Their results are not strange for that if  we suppose $\phi_k^{(n)}(i)=i+\psi_k^{(k)}$ and $\{\psi_k^{(n)}\}$ is a sequence of non-negative mutually independent random variables, \eqref{1.2} can be viewed as a GW process with immigration. From a classical result in \cite{Ka71}, under mild conditions, $\{Y_k(t)\}$ converges to a \textit{continuous state branching process with immigration}. Inspired by this idea, we assume in this paper that
 \beqlb\label{1.3}
 \phi_k^{(n)}(i)=i+\psi_k^{(n)}(i),\ \quad n\geq 0, i\geq 1,
 \eeqlb
 where $\psi_k^{(n)}(i)$ takes non-negative integer. Different from the GW process with immigration, its immigration depends on the current state. Hence, there is a natural conjecture that the limit process is a \textit{continuous state branching process with dependent immigration} (CBDI process).

Let $(\Omega, \mcr{F}, \mcr{F}_t, \mbf{P})$ be a filtered probability space satisfying the usual hypotheses. Let $\{B_t\}$ be an $(\mcr{F}_t)$-Brownian motion.
Let $N_0(\d s,\d z,\d u)$ and $N_1(\d s,\d z,\d u)$ be $\mcr{F}_t$-Poisson random measures with intensities $\d s m(\d z)\d u$ and $\d s\pi(\d z)\d u$ on $(0,\infty)^3$, respectively, where $ (z\wedge z^2)m(\d z)$ is a finite measure and $\pi(dz)$ is a $\sigma$-finite measure. Denote the compensated measures of $\{N_0(\d s,\d z,\d u)\}$  by $\{\tilde{N}_0(\d s,\d z,\d u)\}$. Let $Y_0$ be a non-negative $\mcr{F}_0$-measurable random variable satisfying $\mbf{E}Y_0<\infty$. Let $Y_0$, $\{B_t\}$, $\{N_0(\d s,\d z,\d u)\}$ and $\{N_1(ds,dz,du)\}$ be mutually independent. A CBDI process  $\{Y_t;~t\geq 0\}$ is a non-negative solution to the stochastic integral equation as follows:
\beqlb\label{1.4}
Y_t
\ar=\ar
Y_0+\int_0^t\sqrt{2cY_s}\d B_s+\int_0^t\int_0^{\infty}\int_0^{Y_{s-}} z \tilde{N}_0(\d s,\d z,\d u)\cr\cr
\ar\ar\quad
+\int_0^t (\beta(Y_s)-b Y_s)\d s+\int_0^t\int_0^{\infty}\int_0^{q(Y_{s-},z)}z N_1(\d s,\d z,\d u), \ \quad t\geq 0.
\eeqlb
where $c\geq0$, $b$ are constants, and $x\mapsto\beta(x)$ is a Borel function on $\mbb{R}_+:=[0,\infty)$, and $(x,z)\mapsto q(x,z)$ is a Borel function on $[0,\infty)\times(0,\infty)$.
Moreover,  $\beta(x),q(x,z)$ take non-negative values.
Here and in the sequel, we make the conventions
 \beqnn
\int_a^b = \int_{(a,b]}
 \quad\mbox{and}\quad
\int_a^\infty = \int_{(a,\infty)}
 \eeqnn
for any\ $b\ge a\ge 0$. Let $C_c^{2}(\mbb{R}_{+})$ be the set of bounded continuous real functions on $\mbb{R}_{+}$ with compact support. By It\^o's formula, the generator $L$ of \eqref{1.4} is defined by
\beqlb\label{1.5}
L f(x)
&=&
\big(-bx+\beta(x)\big) f'(x)+x\int_{0}^{\infty}[f(x+z)-f(x)-z f'(x)]
m(dz)\cr\cr
&&\quad+\int_0^{\infty}[f(x+z)-f(x)]q(x,z)\pi(\d z),\ \quad f\in C_c^{2}(\mbb{R}_{+}).
\eeqlb
Throughout this paper, the following assumptions are adopted:
\bitemize
\itm there is a constant $K\geq 0$,
such that
\beqlb\label{1.6}
|\beta(x)|+\int_0^{\infty}q(x,z)z\pi(\d z)\leq K(1+x),\ \quad x\geq 0;
\eeqlb
 \itm
  there is a non-decreasing and concave function
$r:[0,\infty) \mapsto [0,\infty)$ such that $
\int_{0+}r(z)^{-1}\d z=\infty$
 and
 \beqnn
 |\beta(x)-\beta(y)|+\int_0^{\infty}|q(x,z)-q(y,z)|z\pi(\d z)\leq r(|x-y|),\ \quad x,y\geq 0.
 \eeqnn
\eitemize
Under these assumptions, there exists a pathwise unique positive strong solution to (\ref{1.4}) by \cite[Theorem 5.1]{FuL10}.

The paper is organized as follows. In Section 2,  we make some preparations and give some mild conditions which will be used in main theorems. The conditions here can be achieved through constructing probability generating functions. In section 3, based on the martingale problem approach, some estimates are given and lead to the tightness of rescaled processes. In section 4, the conjecture to be a CBDI process is finally verified by using Skorohod  Representative Theorem.


\section{Preliminaries}

\setcounter{equation}{0}
 Let $E_k :=\{ 0, k^{-1}, 2k^{-1},\cdots \}$. Let $g_k$ be the probability generating function of $\xi_{n,i}^{(k)}$. In view of \eqref{1.3}, define the probability generating function of $\psi_k^{(n)}(i)$ by $h_k^{(i)}(s)$. For $0\leq\lambda\leq k$, set
\beqlb\label{2.1}
 R_k(\lambda)=k\gamma_k[g_k(1-\lambda/k)-(1-\lambda/k)],
 \eeqlb
 \beqlb\label{2.2}
F_k(\lambda,x)=\gamma_k[h_k^{(\lfloor kx\rfloor)}(1-\lambda/k)-(1-\lambda/k)].
\eeqlb

We consider conditions as follows:
\bitemize

\itm[{\rm(A)}] The sequence $\{R_k(\lambda)\}$ is uniformly Lipschitz on
each bounded interval and converges to a continuous function $R(\lambda)$
as $k\rightarrow\infty$, where
 \beqlb\label{2.3}
 R(\lambda)= b \lambda +c \lambda^2+\int_0^{\infty}(\e^{-\lambda z}-1+\lambda z)m(dz).
 \eeqlb
\itm[{\rm(B)}] For each $a_1,a_2>0$, the sequence $\{F_k(\lambda,x)\}$
converges to $F(\lambda,x)$ uniformly on $[0,a_1]\times[0,a_2]$ as
$k\rightarrow\infty $, where
\beqlb\label{2.4}
F(\lambda,x)= -\beta(x)\lambda+\int_0^{\infty}(\e^{-\lambda z}-1)q(x,z)\pi(\d z).
\eeqlb
\itm[{\rm(C)}] There exists some positive constant $K_1$ such that
\beqnn
\Big|\frac{\partial F_k}{\partial \lambda}(0,x)\Big|=\frac{\gamma_k}{k}\frac{\d}{\d s}h_k^{(\lfloor kx\rfloor)}(1-)
\leq K_1(1+x), \ \quad x\geq 0, k\geq 1.
\eeqnn
\eitemize
It follows from (C) that there exists some positive constant $K_2$ such that
\beqlb\label{2.6}
|\gamma_k(1-g_k'(1-))|\leq K_2.
\eeqlb
In fact, condition (A) is from \cite[Condition 2.4]{Li20}, and condition (B) is a generalized form of \cite[Condition 5.3]{Li20} in the setting of dependent immigration. Condition (C) is about the first moment, which is viewed as a generalization of \eqref{2.6}.

\btheorem\label{p2.1}
For any function $(R,F)$ with representations (\ref{2.3}) and (\ref{2.4}), respectively, there are  sequences $\{\gamma_k\}$ and $\{(R_k,F_k)\}$ as \eqref{2.1} and \eqref{2.2} satisfying (A-C).
\etheorem
\pf The similar proof was discussed in \cite[Prop 2.6]{Li20}, which should be improved on the present proof. By checking the proof there,  it's sufficient to construct $\{F_k\}$ satisfying (B-C) and  $\gamma_k$ satisfying
\beqlb\label{2.7}
\gamma_k\geq |b|+2ck+\int_0^{\infty}u(1-\e^{-k u})m(\d u).
\eeqlb
Following this step, let $\tilde{F}(\lambda,x)=F(\lambda,x)+\beta(x)\lambda$, and $\tilde{F}_k(\lambda,x)=\tilde{F}(\lambda,\lfloor kx\rfloor/k)\mbf{1}_{[0,k]^2}(\lambda,x)$,
which implies that $\tilde{F}_k(\lambda,x)$ converges to $\tilde{F}(\lambda,x)$ uniformly on the interval $[0,a_1]\times[0,a_2]$ for each $a_1,a_2>0$.
 Next, we need only to adjust $\tilde{\gamma}_k$ to satisfy that for $k\geq 1$,
\beqnn
\tilde{h}_k^{(kx)}(s)=s-\tilde{\gamma}_k^{-1}\tilde{F}\big(k(1-s), x\big),
\ \quad x\in [0,k]\cap E_k
\eeqnn
is a probability generating function.
Above all, $s\mapsto \tilde{h}_k^{(kx)}(s)$ is an analytic
function. By elementary calculations, to ensure that
\beqnn
\frac{\d^n}{\d z^n}\tilde{h}_k^{(kx)}(0)\geq 0,\ \quad n\geq 0,
\eeqnn
 it's sufficient to show that
\beqnn
\tilde{\gamma}_k\geq  \int_0^{\infty}k z\e^{-kz}q(x,z)m(\d z),
\ \quad x\leq k.
\eeqnn
Consequently, (\ref{1.6}) implies that there is a sequence $\{\tilde{\gamma}_k\}$ such that
$\tilde{h}_k^{(kx)}(s)$ is a probability generating function for
$x\in E_k\cap[0,k],k\geq 1$. On the other hand, since \eqref{1.6} implies
$\beta(x)\leq \hat{K}\sqrt{k}$ for $x\leq \sqrt{k}$ and some positive constant $\hat{K}$
independent of $k$,  we define for $x\leq\sqrt{k}$,
\beqnn
\hat{h}_k^{(kx)}(s)=
 1-\frac{2\beta(\lfloor kx\rfloor/k)}{3\sqrt{k}\hat{K}}
 +\frac{\beta(\lfloor kx\rfloor/k)}{3\sqrt{k}\hat{K}}s+\frac{\beta(\lfloor kx\rfloor/k)}{3\sqrt{k}\hat{K}}s^2;
 \eeqnn
and for $x>\sqrt{k}$, let $\hat{h}_k^{(kx)}(s)=1$. For arbitrary $k\geq 1$, define constant $\gamma_k=3k^{3/2}\hat{K}$ and
\beqnn
\hat{F}_k(\lambda, x)\ar:=\ar
\hat{\gamma}_k[\hat{h}_k^{(\lfloor kx\rfloor)}(1-\lambda/k)-1]
\mbf{1}_{[0,k]\times[0,\sqrt{k}]}(\lambda,x)\cr\cr
\ar=\ar
-\big[\beta(\lfloor kx\rfloor/k)\lambda-\frac{1}{3k}\beta(\lfloor kx\rfloor/k)\lambda^2\big]
\mbf{1}_{[0,k]\times[0,\sqrt{k}]}.
\eeqnn
It's easy to verify that $\hat{F}_k^{(kx)}(s)$ tends to $-\beta(x)\lambda$ uniformly on $[0,a_1]\times[0,a_2]$ as
$k\rightarrow\infty $. In the end, set
\beqlb\label{2.8}
\gamma_k=\tilde{\gamma}_k+\hat{\gamma}_k,\ \quad h^{(kx)}_k(s)
=\gamma_k^{-1}[\tilde{\gamma}_k\tilde{h}_k^{(kx)}(s)+\hat{\gamma}_k\hat{h}_k^{(kx)}(s)].
\eeqlb
Consequently, for $k\geq 1$,
\beqnn
F_k(\lambda,x)\ar=\ar\tilde{F}_k(\lambda,x)+\hat{F}_k(\lambda,x), \cr\cr
\ar=\ar-\big[\beta(\frac{\lfloor kx\rfloor}{k})\lambda-\frac{1}{3k}\beta(\frac{\lfloor kx\rfloor}{k})\lambda^2\big]
\mbf{1}_{[0,k]\times[0,\sqrt{k}]}\cr\cr
\ar\ar\qquad+\int_0^{\infty}(\e^{-\lambda z}-1)q(x,z)\pi(\d z)\mbf{1}_{[0,k]^2}
\eeqnn
satisfies (B). By elementary calculations and combining with \eqref{1.6}, we can verify  (C). Finally, from \eqref{2.8}, $\gamma_k$ also satisfies \eqref{2.7}, which ensures (A). Then we get the desired result.
\qed

For convenience, set $e_\lambda(x):=\e^{-\lambda x}$ for $\lambda\geq 0$ and $x\geq 0$.
Let $D_0$ be the linear hull of $\{e_{\lambda}(x):\lambda\geq 0\}$. Denote the set of
bounded measurable real functions on $\mbb{R}_{+}$ by $b(\mbb{R}_{+})$.  Next, we introduce an analytical conclusion
for the sake of clarity.

\blemma\label{l2.2}
For $f\in C_c^2(\mbb{R}^{+})$, there exists a sequence of functions $f_n$
in $D_0$, such that
\beqnn
f_n\rightarrow f,\qquad f'_n\rightarrow f',\qquad f''_n\rightarrow f''
\eeqnn
 uniformly
on $\mbb{R}_{+}$, as $n\rightarrow\infty$.
\elemma
\pf
For $f\in C_c^2(\mbb{R}^{+})$, define a function $p$ on $[0,1]$ by
 \beqnn
p(x)=
 \left\{
\begin{array}{ll}
 f(-\log(x)),&x>0, \cr
 0,\ &x=0.
\end{array}
\right.
 \eeqnn
By simple calculations, $p\in C^2[0,1]$. For a real function $p$ on
$[0,1]$, its Bernstein polynomial is given by
\beqnn
B_n(p,x)=\sum_{r=0}^n p(r/n)\binom{n}{r} x^r(1-x)^{n-r}.
\eeqnn
Let
\beqnn
f_n(x)=\sum_{k=1}^n \binom{n}{k} f\big(-\log(k/n)\big)\e^{-kx}(1-\e^{-x})^{n-k},
\eeqnn
which implies that $f_n\in D_0$ and
\beqlb\label{2.9}
f_n(x)=B_n(p,\e^{-x}).
\eeqlb
By taking derivatives on both sides of (\ref{2.9}), we have
\beqnn
f_n'(x)=-\e^{-x}B_n'(p,\e^{-x}).
\eeqnn
In fact, as a result of \cite[Theorem 7.16]{Ph03}, $\lim_{n\rightarrow\infty}f_n=f$
uniformly on $\mbb{R}_{+}$. Besides, $B'_n(p,\e^{-x})$ converges to $p'(\e^{-x})$
uniformly on $\mbb{R}_{+}$ as $n\rightarrow\infty$. Hence, $f_n'(x)\rightarrow f'(x)$
uniformly on $\mbb{R}_{+}$ as $n\rightarrow\infty$. The same argument leads to the
desired result.
\qed


\section{Discrete martingale  and tightness}

\setcounter{equation}{0}
In this section, we construct the discrete martingale and prove the tightness
of the rescaled processes.

Let $D([0,\infty),\mbb{R}_{+})$ be the space of c\`adl\`ag functions
$\omega: [0,\infty)\mapsto[0,\infty)$. For $\lambda>0$, the distance
$\rho_{\lambda}$ is defined by
  \beqnn
  \rho_{\lambda}(x,y)=|\e^{-\lambda x}-\e^{-\lambda y}|,\ \quad x,y\in[0,\infty)
  \eeqnn
Denote the one-step transition matrix of $Z_k(n)/k$  to be $T_k$. Suppose that the process $\{Z_k(n)/k\}_n$ is adapted to a filtration $(\mathscr{G}_n)_{n\geq 1}$ for each $k$. Define $A_k$ as the discrete generator of $\{Y_k(t)\}$; for more details on discrete generators, see \cite[pp 230-233]{Et86}. Then
\beqnn
A_k f=\gamma_k(T_k-I)f, \ \quad f\in b(\mbb{R}_+),
\eeqnn
where $I(x)=x$ for $x\geq 0$. Based on discrete generators, we can construct the discrete martingale problem.

\blemma\label{l3.1}
For $f\in b(\mbb{R}_+)$, set
\beqnn
 M_k^{f}(n):=\gamma_k f\big(\frac{Z_k(n)}{k}\big)-\gamma_k f\big(\frac{Z_k(0)}{k}\big)
 -\sum_{i=0}^{n-1}A_k f\big(\frac{Z_k(i)}{k}\big).
 \eeqnn
Then $\{M_k^{f}(n)\}_n$ is a $(\mathscr{G}_n)$-martingale.
\elemma
\pf
It follows directly from Markov property.\qed

 It is obvious by elementary calculations that
\beqnn
A_k e_{\lambda}(x)=\gamma_k \mbf{E}
[g_k(\e^{-\lambda /k})^{\phi_k(kx)}-\e^{-\lambda x}].
\eeqnn
Recall from \eqref{1.5}. Then for $\lambda\geq0$,
\beqlb\label{3.1}
L e_{\lambda}(x)=x\e^{-\lambda x}R(\lambda)+\e^{-\lambda x}F(\lambda,x),
\ \quad x\geq 0.
\eeqlb
By (\ref{1.6}), $Le_{\lambda}(x)\rightarrow 0$ as $x\rightarrow\infty$.
Hence $Ae_{\lambda}(x)$ is bounded on $\mbb{R}_{+}$. Based on assumptions on probability generating functions, we have the following estimate.

\btheorem\label{t3.2}
Suppose that (A-C) hold. Then for $\lambda\geq 0$, we have
\beqnn
\lim_{k\rightarrow\infty}\sup_{x\in E_k}|L e_{\lambda}(x)-A_k e_{\lambda}(x)|=0.
\eeqnn
\etheorem
\pf
Observe that
\beqnn
A_k e_{\lambda}(x)=\gamma_k\mbf{E}[g_k(\e^{-\lambda /k})^{\phi_k(kx)}
-g_k^{kx}(\e^{-\lambda /k})]+\gamma_k[g_k^{kx}(\e^{-\lambda /k})-\e^{-\lambda x}].
\eeqnn
Let
\beqlb\label{3.2}
B_k(\lambda, x)\ar=\ar
\gamma_k\mbf{E}[g_k(\e^{-\lambda /k})^{\phi_k(kx)}-g_k^{kx}(\e^{-\lambda /k})]\cr\cr
\ar=\ar \gamma_k g_k^{kx}(\e^{-\lambda /k})[h_k^{(kx)}(g_k(\e^{-\lambda/k}))-1],
\eeqlb and
\beqnn
C_k(\lambda, x)=\gamma_k[g_k^{kx}(\e^{-\lambda /k})-\e^{-\lambda x}].
\eeqnn
Similar to \cite[Theorem 2.1]{Li06}, there is a more precise approximation:
\beqlb\label{3.3}
\lim_{k\rightarrow \infty}\sup_{x\in E_k}\e^{\lambda_0 x}
|C_k(\lambda, x)-x \e^{-\lambda x}R(\lambda)|=0, \ \quad\lambda_0<\lambda.
\eeqlb
Set
\beqnn
 u_k(\lambda):=k[1-g_k(\e^{-\lambda/k})].
\eeqnn
It's easy to check that for every $a\geq 0$, $u_k(\lambda)\rightarrow\lambda$
uniformly on $[0,a]$ as $k\rightarrow\infty$. As a result,
$(1-u_k(\lambda)/k)^{kx}$ converges to $\e^{-\lambda x}$ uniformly on
$\mbb{R}_{+}$ as $k\rightarrow\infty$.
Consequently, combining it with (\ref{3.2}),
\beqnn
B_k(\lambda, x)
=(1-u_k(\lambda)/k)^{kx} F_k\big(u_k(\lambda),x\big),
\eeqnn
which follows from (2.B) that $B_k(\lambda,x)$
converges to $\e^{-\lambda x}F(\lambda,x)$ uniformly on $[0,M]$
as $k\rightarrow\infty$, for $M>0$.
On the other hand, (C) yields
\beqlb\label{e2.5}
|F_k(\lambda,x)|\leq K_1\lambda(1+x), \quad\lambda,x\geq 0,k\geq 1.
\eeqlb
Therefore, there is a constant $\tilde{C}>0$, such that
\beqnn
\sup_{k\geq 1}|F(\lambda,x)-F_k\big(u_k(\lambda),x\big)|\leq \tilde{C} \e^{\lambda_0 x}.
\eeqnn
Define
\beqnn
\epsilon_k(x):=e^{\lambda_0 x}[C_k(\lambda, x)-x \e^{-\lambda x}R(\lambda)].
\eeqnn
By (\ref{3.3}), $\epsilon_k(x)$ converges to $0$ uniformly, as
$k\rightarrow \infty$. Therefore,
\beqnn
\ar\ar
\limsup_{k\rightarrow\infty}\sup_{x\in E_k\cap[M,\infty]}
g^{kx}_k(\e^{-\lambda/k})
|F(\lambda,x)-F_k\big(u_k(\lambda),x\big)|\cr\cr
\ar\ar\qquad
\leq\limsup_{k\rightarrow\infty}\sup_{x\in E_k\cap[M,\infty]}
[\gamma_k^{-1}x\e^{-\lambda x}R(\lambda)
+\gamma_k^{-1}\e^{-\lambda_0 x}\epsilon_k(x)+\e^{-\lambda x}]\tilde{C}\e^{\lambda_0 x}
\cr\cr
\ar\ar\qquad
\leq \tilde{C}\e^{(\lambda_0-\lambda)M},
\eeqnn
which yields
\beqnn
\lim_{M\rightarrow\infty}\varlimsup_{k\rightarrow\infty}
\sup_{x\in E_k\cap[M,\infty]}
|B_k(\lambda, x)-F(\lambda,x)\e^{-\lambda x}|=0.
\eeqnn
That gives the desired result along with the result that $B_k(\lambda,x)$
converges to $\e^{-\lambda x}F(\lambda,x)$ uniformly on $[0,M]$
as $k\rightarrow\infty$, for $M>0$.
\qed

\bremark
Since we can't  prove that $D_0$ is  a core for $L$, the result for weak convergence couldn't be directly obtained from \cite[Corollary 8.9]{Et86} similar to \cite[Theorem 2.1]{Li06}. On the other hand, the result of Gonz\'alez and del Puerto \cite{Go12} is also from another result in  \cite[Corollary 8.9]{Et86}. However, it's not adapted in our scene, for  our estimate is restricted on $f=e_{\lambda}$. In the following pages, we will prove tightness and use Skorohod Representative Theorem to avoid the barrier.
\eremark

In the rest of this section, we aim at the proof of tightness of $\{Y_k(t)\}$. For convenience,
we introduce some notations before that. For a fixed constant $T>0$,
we consider a sequence of stopping times $\tau_k$. Let $\delta_k$ be a sequence of positive
constants that tends to $0$ as $k$ tends to infinity. Suppose
\beqnn
0\leq\tau_k<\tau_k+\delta_k\leq T.
\eeqnn
\blemma\label{l3.3}
Suppose that (A-C) hold. Then
\beqnn
\lim_{k\rightarrow\infty}\mbf{E}[\rho_{\lambda}^2\big(Y_k(\tau_k+\delta_k),Y_k(\tau_k)\big)]=0.
\eeqnn
\elemma
\pf
By Lemma \ref{l3.1}, for $f\in b(\mbb{R}_{+})$,
\beqlb\label{3.4}
\mbf{E}[f\big(Y_k(t)\big)]
\ar=\ar
\mbf{E}\Big[f\big(k^{-1}Z_k(\lfloor \gamma_k t \rfloor)\big)\Big]\cr
\ar=\ar
\mbf{E}\Big[f\big(Y_k(0)\big)
+\sum_{i=0}^{\lfloor \gamma_k t \rfloor-1}\gamma_k^{-1}
A_k f\big(k^{-1}Z_k(i)\big)\Big]\cr
\ar=\ar
\mbf{E}\big[f\big(Y_k(0)\big)\big]+\mbf{E}
\Big[\int_0^{\lfloor \gamma_k t \rfloor/\gamma_k}A_k f\big(Y_k(s)\big)\d s\Big].
\eeqlb
For $\lambda>0$,
\beqlb\label{3.5}
\ar\ar\mbf{E}\big[|\e^{-\lambda Y_k(\tau_k+\delta_k)}-\e^{-\lambda Y_k(\tau_k)}|^2\big]
\cr\cr
\ar\ar\qquad=
\mbf{E}[e_{2\lambda}\big(Y_k(\tau_k)\big)-2e_{\lambda}\big(Y_k(\tau_k)+Y_k(\tau_k+\delta_k)\big)
+e_{2\lambda} \big(Y_k(\tau_k+\delta_k)\big)]
\cr\cr
\ar\ar\qquad\leq
I_1+I_2+I_3,
\eeqlb
where
\beqnn
I_1
\ar:=\ar
\big|\mbf{E}\big[\e^{-2\lambda Y_k(\tau_k+\delta_k)}-\e^{-2\lambda Y_k(\tau_k)}\big]\big|;\cr\cr
I_2
\ar:=\ar
\Big|\mbf{E}\Big[2\e^{-\lambda Y_k(\tau_k)}\int^{{\lfloor \gamma_k (\tau_k+\delta_k) \rfloor}/
{\gamma_k}}_{{\lfloor \gamma_k \tau_k \rfloor}/{\gamma_k}}
A_k \e^{-\lambda Y_k(u)}\d u\Big]\Big|;\cr\cr
I_3
\ar:=\ar
\big|\mathbf{E}\big[2\gamma_k^{-1}\e^{-\lambda Y_k(\tau_k)}
[M_k^{e_{\lambda}(\cdot)}(\lfloor \gamma_k (\tau_k+\delta_k) \rfloor)
-M_k^{e_{\lambda}(\cdot)}(\lfloor \gamma_k \tau_k \rfloor)]\big]\big|.
\eeqnn
Then a simple application of  (\ref{3.4}) yields
\beqlb\label{3.6}
I_1
\ar=\ar
\Big|\mbf{E}\Big[\int^{\lfloor \gamma_k (\tau_k+\delta_k) \rfloor/\gamma_k}
_{\lfloor \gamma_k \tau_k \rfloor/\gamma_k}A_k e_{2\lambda} \big(Y_k(s)\big)\d s\Big]\Big|\cr\cr
\ar\leq\ar
\mbf{E}\Big[\int^{\lfloor \gamma_k (\tau_k+\delta_k) \rfloor/\gamma_k}
_{\lfloor \gamma_k \tau_k \rfloor/\gamma_k}
|A_k e_{2\lambda}\big(Y_k(s)\big)-A e_{2\lambda}\big(Y_k(s)\big)|\d s\Big]\cr\cr
\ar\ar\qquad
+\mbf{E}\Big[\int^{\lfloor \gamma_k (\tau_k+\delta_k) \rfloor/\gamma_k}
_{\lfloor \gamma_k \tau_k \rfloor/\gamma_k}\big|A e_{2\lambda}\big(Y_k(s)\big)\big|\d s\Big]\cr\cr
\ar\leq\ar C_1 \delta_k,
\eeqlb
where $C_1$ is a positive constant.
The same argument implies that
\beqlb\label{3.7}
I_2
\ar\leq\ar
\Big|2\mbf{E}\Big[\int^{\lfloor \gamma_k (\tau_k+\delta_k) \rfloor/\gamma_k}
_{\lfloor \gamma_k \tau_k \rfloor/\gamma_k}A_k e_{\lambda} \big(Y_k(s)\big)\d s\Big]\Big|\cr\cr
\ar\leq\ar C_2 \delta_k,
\eeqlb
for a positive constant $C_2$.
On the other hand,
\beqnn
\Omega
\ar=\ar\big\{\lfloor \gamma_k (\tau_k+\delta_k) \rfloor
=\lfloor \gamma_k\tau_k\rfloor+\lfloor\gamma_k\delta_k\rfloor\big\}
\bigcup
 \big\{\lfloor \gamma_k (\tau_k+\delta_k) \rfloor
 =\lfloor \gamma_k\tau_k\rfloor+\lfloor\gamma_k\delta_k\rfloor+1\big\}\cr\cr
 \ar:=\ar \Omega_1+\Omega_2.
\eeqnn
Both $\Omega_1$ and $\Omega_2$ are $\mcr{G}_{\lfloor \gamma_k\tau_k\rfloor}$-measurable.
Observe that $Y_k(\tau_k)$ is also $\mcr{G}_{\lfloor \gamma_k\tau_k\rfloor}$-measurable.
Then, from the results of Lemma \ref{l3.1} and Doob's Stopping Theorem, it follows that
\beqnn
\mbf{E}\big\{\e^{-\lambda Y_k(\tau_k)}
\big[M_k^{e_{\lambda}(\cdot)}\big(\lfloor \gamma_k (\tau_k+\delta_k) \rfloor\big)
-M_k^{e_{\lambda}(\cdot)}(\lfloor \gamma_k \tau_k \rfloor)\big];\Omega_i\big\}=0,\ \quad i=1,2.
\eeqnn
Thus, $I_3=0$. Together with (\ref{3.5}), (\ref{3.6}) and (\ref{3.7}),
since $\lim_{k\rightarrow\infty}\delta_k=0$, we obtain the desired result.
\qed

In the following, we need the first moment condition for $Y_k(0)$.
\bitemize
 \itm[{\rm(D)}] $\sup_{k\geq 1}\mbf{E}[Y_k(0)]<\infty$.
 \eitemize
 Based on it, a precise estimate on $\mbf{E}[Y_k(\tau_k)]$ is obtained.
\blemma\label{l3.4}
Suppose that (A), (C) and (D) hold. Then there exists some constant $K_3\geq 0$, such that
\beqnn
\mbf{E}[Y_k(\tau_k)]\ \mbox{\text{and}}\ \mbf{E}[Y_k(\tau_k+\delta_k)]\leq K_3,\ \quad k\geq 1.
\eeqnn
\elemma
\pf
\textbf{Step 1.} We prove that there exists some constant $K_4\geq 0$, such that
\beqlb\label{3.8}
\mbf{E}[Y_k(s)]\leq K_4,\ \quad s\leq T, k\geq 1.
\eeqlb
In fact, by (C) and (\ref{2.6}),
\beqnn
\mbf{E}[Z_k(n+1)/k]\ar=\ar g_k'(1-)\mbf{E}[\phi_k(Z_k(n))/k]\cr\cr
\ar=\ar g_k'(1-)\mbf{E}[\frac{1}{k}\big(\frac{\d}{\d s}h_k^{(Z_k(n))}(1-)+Z_k(n)\big)]\cr\cr
\ar\leq\ar g_k'(1-)\big(\mbf{E}[Z_k(n)/k](1+K_1/\gamma_k)+K_1/\gamma_k\big)\cr\cr
\ar\leq\ar (1+K_2/\gamma_k)\big(\mbf{E}[Z_k(n)/k](1+K_1/\gamma_k)+K_1/\gamma_k\big).
\eeqnn
 By induction,
\beqlb\label{3.9}
\mbf{E}[Y_k(t)]\ar\leq\ar\{(1+K_2/\gamma_k)(1+K_1/\gamma_k)\}
^{\lfloor \gamma_k t\rfloor}\mbf{E}[Y_k(0)]\cr\cr
\ \quad \ar\ar+K_1/\gamma_k\frac{\{(1+K_2/\gamma_k)(1+K_1/\gamma_k)\}
^{\lfloor \gamma_k t\rfloor}-1}{(1+K_2/\gamma_k)(1+K_1/\gamma_k)-1}
\eeqlb
A simple calculation yields (\ref{3.8}).

\textbf{Step 2.}  By (C) and (\ref{2.6}),
\beqnn
|A_k I(x)|\ar=\ar
\Big|\gamma_k\mbf{E}\Big[\frac{1}{k}\sum_{i=1}^{\phi_k(kx)}\xi_{n,i}^{(k)}-x\Big]\Big|\cr\cr
\ar=\ar
|k^{-1}\gamma_kg_k'(1-)\frac{\d}{\d s}h_k^{(kx)}(1-)+x\gamma_k(g_k'(1-)-1)|\cr\cr
\ar\leq\ar
|k^{-1}\gamma_kg_k'(1-)\frac{\d}{\d s}h_k^{(kx)}(1-)|+|x\gamma_k\big(1-g_k'(1-)\big)|\cr\cr
\ar\leq\ar
(1+\frac{K_2}{\gamma_k})K_1(1+x)+x K_2\cr\cr
\ar\leq\ar K_5 (1+x), \ \quad k\geq 1,~x\geq0,
\eeqnn
for some positive constant $K_5$,
which follows from \eqref{3.9} that $A_k I(Y_k(s))$ is integrable for $s\leq t,k\geq 1$.
The same argument as the proof in Lemma \ref{l3.1} for $f=I$ leads $M^{I}_k(n)$ to be a martingale.
By Doob's Stopping Theorem,
\beqnn
\mbf{E}[Y_k(\tau_k)]\ar=\ar
\mbf{E}[Y_k(0)]+\mbf{E}[\int_0^{\tau_k}A_k Id(Y_s) \d s]\cr\cr
\ar\leq\ar
\sup_{k\geq 1}\mbf{E}[Y_k(0)]+\int_0^T K_5\mbf{E}[1+Y_k(s)]\d s\cr\cr
\ar\leq\ar
\sup_{k\geq 1}\mbf{E}[Y_k(0)]+T K_5 (1+K_4).
\eeqnn
The same argument leads to $\mbf{E}[Y_k(\tau_k+\delta_k)]\leq K_3$, 
which completes the proof.
\qed

Now, we give the tightness using Aldous' criterion.
\btheorem\label{t3.5}
Suppose that (A-D) hold. Then the process
$\{Y_k(t): t\geq 0\}_{k\geq 1}$ is tight in $D([0,\infty),\mbb{R}_{+})$.
\etheorem
\pf
Firstly, it follows from (\ref{3.8}) that $\{Y_k(t)\}$ is tight for a fix $t\geq 0$. Next,
 for a fixed constant $M>0$, for $\epsilon>0$, $|a-b|>\epsilon$ and  $0\leq a,b\leq M$, 
 \beqnn
|\e^{-\lambda a}-\e^{-\lambda b}|\geq \lambda \e^{-\lambda M}\epsilon.
\eeqnn
Hence, by a simple calculation,
\beqnn
 \ar\ar\mbf{P}\{|Y_k(\tau_k+\delta_k)-Y_k(\tau_k)|>\epsilon;~Y_k(\tau_k)\vee Y_k(\tau_k+\delta_k)\leq M\}\cr\cr
&&\qquad\qquad\qquad\leq
(\lambda\epsilon)^{-2}\e^{2\lambda M}\mbf{E}[|\e^{-\lambda Y_k(\tau_k+\delta_k)}-\e^{-\lambda Y_k(\tau_k)}|^2].
\eeqnn
By Lemma \ref{3.3},
\beqnn
\lim_{k\rightarrow\infty}\mbf{P}\{|Y_k(\tau_k+\delta_k)-Y_k(\tau_k)|>\epsilon;~ Y_k(\tau_k)\vee Y_k(\tau_k+\delta_k)\leq M\}=0.
\eeqnn
On the other hand, by Lemma \ref{3.4},  
\beqnn
\mbf{P}(Y_k(\tau_k)\geq M)\leq K_3/M\ \mbox{\text{and}}\ \mbf{P}(Y_k(\tau_k+\delta_k)\geq M)\leq K_3/M
\eeqnn
As a result,
\beqnn
\ar\ar\mbf{P}\{|Y_k(\tau_k+\delta_k)-Y_k(\tau_k)|>\epsilon\}\cr\cr
\ar\ar\qquad\qquad\leq
\mbf{P}\{|Y_k(\tau_k+\delta_k)-Y_k(\tau_k)|>\epsilon; Y_k(\tau_k)\vee Y_k(\tau_k+\delta_k)\leq M\}
\cr\cr\ar\ar\qquad\qquad\qquad\qquad\qquad+\mbf{P}(Y_k(\tau_k)\geq M)+\mbf{P}(Y_k(\tau_k+\delta_k)\geq M)\cr\cr
\ar\ar\qquad\qquad\leq
\mbf{P}\{|Y_k(\tau_k+\delta_k)-Y_k(\tau_k)|>\epsilon; Y_k(\tau_k)\vee Y_k(\tau_k+\delta_k)\leq M\}+2K_3/M.
\eeqnn
Let $k\rightarrow\infty$ and $M\rightarrow\infty$ following, then we obtain
\beqnn
\lim_{k\rightarrow\infty}\mbf{P}\{|Y_k(\tau_k+\delta_k)-Y_k(\tau_k)|>\epsilon\}=0.
\eeqnn
Finally, the tightness of
$\{Y_k(t): t\geq 0\}_{k\geq 1}$ in $D([0,\infty),\mbb{R}_+)$ follows from
Aldous' criterion in \cite[Theorem 1]{Al78}.
\qed

\section{Weak convergence}

\setcounter{equation}{0}
In this section, we build the relations between the weak solution of (\ref{1.4})
and its corresponding martingale problem. By an application of Skorokhod Representative Theorem,
the weak limit of rescaled processes is proved to be a weak solution of (\ref{1.4}).

\btheorem\label{p4.1}
A positive c\`adl\`ag process $\{Y_t:t\geq 0\}$ is a weak solution of (\ref{1.4})
with initial value $Y_0$ if and only if for every $f\in C_c^{2}(\mbb{R}_{+})$ ,
\beqlb\label{4.2}
f(Y_t)=f(x_0)+\int_0^t L f(Y_s)\d s + \text{local mart}, \ \quad y\geq 0.
\eeqlb
\etheorem
\pf
The proof is a modification of that in \cite[Proposition 4.2]{FuL10} which needs a stronger condition and that in \cite[Theorem 5.1]{Li19} which consider a different domain of generator. For rigorousness, we give a brief proof for the different part from their proofs. Suppose that (\ref{4.2}) holds for
every $f\in C_c^{2}(\mbb{R}_{+})$.
We introduce a non-decreasing sequence of functions $f_n\in C_c^{2}(\mbb{R}_{+})$
such that $f_n(x)=x$, for $0\leq x\leq n$ and $f_n'(x)\leq 1$ for $x\geq0$.
Let $\tau_n=\inf\{t>0, Y_t\geq n\}$. It follows from (\ref{4.2}) that
\beqlb\label{4.3}
M_n(t\wedge\tau_m):=f_n(Y_{t\wedge\tau_m})-f_n(Y_0)-\int_0^t Lf_n(Y_{s\wedge\tau_m-})\d s,
\ \quad m\leq n
\eeqlb
is a martingale. Consequently, taking expectations above, we obtain
\beqnn
\mbf{E}[f_n(Y_{t\wedge\tau_m})]
\ar=\ar
\mbf{E}f_n(Y_0)+\int_0^t\mbf{E}[-b Y_{s\wedge\tau_m-}+\beta(Y_{s\wedge\tau_m-})]\d s\cr\cr
\ar\ar\quad
+\int_0^t\int_0^{\infty}\mbf{E}[f_n(Y_{s\wedge\tau_m-}+z)-Y_{s\wedge\tau_m-}-z]
\d sm(\d z)\cr\cr
\ar\ar\quad
+\int_0^t\int_0^{\infty}\mbf{E}[\big(f_n(Y_{s\wedge\tau_m-}+z)-Y_{s\wedge\tau_m-}\big)
q(Y_{s\wedge\tau_m-},z)]\d s\pi(\d z)\cr\cr
\ar\leq\ar
\mbf{E}Y_0+|b|mt+K(1+m)t.
\eeqnn
Hence, by monotone convergence, $Y_{t\wedge\tau_m}$ is integrable for $m\geq1$.
Letting $n\rightarrow\infty$ in (\ref{4.3}),
by monotone convergence, we obtain that
\beqnn
Y_{t\wedge\tau_m}-Y_0-\int_0^t[b Y_{s\wedge\tau_m-}+\beta(Y_{s\wedge\tau_m-})]\d s-\int_0^t \int_0^{\infty}z q(Y_{s\wedge\tau_m-},z)\pi(\d z)\d s
\eeqnn
is a martingale. We omit the rest proof, for it's the same as that
in \cite[Proposition 4.2]{FuL10}.
\qed

Finally, combining all the results above, we can obtain the weak convergence.
\btheorem\label{t4.2}
Suppose that (A-D) hold and $Y_k(0)$ converges in distribution to $Y_0$ as $k\rightarrow\infty$. Then $\{Y_k(t):t\geq 0\}_k$
converges in distribution on $D([0,\infty),\mbb{R}_{+})$ to $\{Y_t:t\geq 0\}$, where
$\{Y_t:t\geq 0\}$ is a weak solution to (\ref{1.4}) with initial value $Y_0$.
\etheorem
\pf
The definition of $M_k(n)$ implies that
\beqlb\label{4.4}
f(Y_k(t))
\ar=\ar
f(Y_k(0))
+\sum_{i=0}^{\lfloor \gamma_kt\rfloor-1}A_k f\big(\frac{Z_k(i)}{k}\big)
+\gamma_k^{-1}M_k^f(\lfloor \gamma_kt\rfloor)\cr\cr
\ar=\ar
f(Y_k(0))
+\int_0^{\lfloor \gamma_kt\rfloor/\gamma_k}A_k f(Y_k(s))\d s
+\gamma_k^{-1}M_k^f(\lfloor \gamma_kt\rfloor),\ \quad f\in b(\mbb{R}_{+}).
\eeqlb
Let $P^{(n)}$ and $P$ be the distributions of $Y_n$ and $Y$ in $D$,
respectively. By Theorem \ref{t3.5}, $\{Y_n\}_n$ is relatively compact.
Then we can find a probability measure $Q$ on $D([0,\infty),\mbb{R}_{+})$ and
a subsequence $P^{(n_i)}$ such that $Q=\lim_{i\rightarrow\infty}P^{(n_i)}$.
By Skorokhod Representative Theorem, there exists a probability space
$(\tilde{\Omega},\tilde{\mathscr{F}},\tilde{\mbf{P}})$ on which are defined
c\`adl\`ag processes $\{X_t:t\geq 0\}$ and $\{X_t^{(n_i)}: t\geq 0\}$
such that
 \bitemize
\itm[{\rm(i)}] the distributions of $X$ and $X^{(n_i)}$ on $D([0,\infty),\mbb{R}_{+})$
are $P$ and $P^{(n_i)}$, respectively;
\itm[{\rm(ii)}] the limit holds
\beqnn
\lim_{i\rightarrow \infty} X^{(n_i)}= X,\ \quad \tilde{\mbf{P}}-a.s.
\eeqnn
\eitemize
Firstly, taking $f=e_{\lambda}$ in (\ref{4.4}), we have
 \beqlb\label{4.5}
e_{\lambda}(X^{(n_i)}_t)=
e_{\lambda}(X^{(n_i)}_0)
+\int_0^{\lfloor \gamma_{n_i}t\rfloor/\gamma_{n_i}}
A_{n_i}e_{\lambda}(X^{(n_i)}_s)\d s
+\gamma_{n_i}^{-1}M_{n_i}^{e_{\lambda}(\cdot)}(\lfloor \gamma_{n_i}t\rfloor).
\eeqlb
Observe that
\beqlb\label{4.6}
\ar\ar\int_0^{\lfloor \gamma_{n_i}t\rfloor/\gamma_{n_i}}
|A_{n_i}e_{\lambda}(X^{(n_i)}_s)-L e_{\lambda}(X_s)|\d s\cr\cr
\ar\ar\quad\leq
\int_0^{\lfloor \gamma_{n_i}t\rfloor/\gamma_{n_i}}
|A_{n_i}e_{\lambda}(X^{(n_i)}_s)-L e_{\lambda}(X^{(n_i)}_s)|\d s\cr\cr
\ar\ar\qquad+\int_0^{\lfloor \gamma_{n_i}t\rfloor/\gamma_{n_i}}
|L e_{\lambda}(X^{(n_i)}_s)-L e_{\lambda}(X_s)|\d s.
\eeqlb
By Theorem \ref{t3.2}, the first integral on the right side of (\ref{4.6})
converges to $0$ as $i\rightarrow\infty$. Let
\beqnn
D(x):=\{t>0:\tilde{\mbf{P}}\{X(t-)=X(t)\}=1\}.
\eeqnn
Then from \cite[p.118]{Et86}, (\rm{ii}) implies that for each $t\in D(x)$,
as $i\rightarrow\infty$, $X^{(n_i)}_t$ converges to $X_t$.
 Observe from \cite[p.131]{Et86} that the set
$\mbb{R}_{+}/D(x)$ is at most countable.
 Therefore, the second integral on the
right side of (\ref{4.6}) also converges to $0$
as $i\rightarrow\infty$. On the other hand, since $Le_{\lambda}(\cdot)$ is
bounded, (\ref{4.6}) also
implies that the second term of the right side hand of
(\ref{4.5}) is uniformly bounded for all $i\geq 1, \omega\in\tilde{\Omega}$.
Consequently, taking limits in (\ref{4.5}) and using
bounded convergence theorem, we have
\beqlb\label{4.7}
M_t^f=f(X_t)-f(X_0)-\int_0^t L f(X_s)\d s, \qquad f=e_{\lambda}
\eeqlb
is a martingale bounded on each bounded time interval.
Next, for $f\in C_c^2(\mbb{R}_{+})$, let $f_m$ be given in Lemma \ref{l2.2}.
In view of (\ref{1.5}), $L f_m(x)$ converges to $L f(x)$ uniformly on each
bounded interval. As a linear span of $\{e_{\lambda}(x)\}$, by (\ref{4.7}),
\beqlb\label{4.8}
f_m(X_t)=f_m(X_0)+\int_0^t L f_m(X_s)\d s+M_t^{f_m}.
\eeqlb
Let $\tilde{\tau}_M:=\inf\{t>0,X_t\geq M\}$.
Since $X_t$ is a c\`adl\`ag process,  $\tilde{\tau}_M\rightarrow\infty$ a.s.,
as $M\rightarrow\infty$.
Replacing $t$ with $t\wedge\tilde{\tau}_M$, and taking limits on both sides of
(\ref{4.8}), we use the same argument to obtain that
 \beqnn
M_{t\wedge\tilde{\tau}_M}^f=f(X_{t\wedge\tilde{\tau}_M})-f(X_0)
-\int_0^t L f(X_{s\wedge\tilde{\tau}_M-})\d s
\eeqnn
is a martingale bounded on each bounded time interval. Hence, $X_t$
solves the martingale problem (\ref{4.2}). By Theorem \ref{p4.1}, $X_t$
is a weak solution of (\ref{1.4}). By \cite[Theorem 137]{Si05}, the pathwise
uniqueness of (\ref{1.4}) implies the uniqueness of distributions. Therefore,
$Q=P$, and $\lim_{i\rightarrow\infty}P^{(n_i)}=P$, which completes the desired
result.
\qed



\end{document}